\documentclass[12pt]{article}

\usepackage{latexsym}
\usepackage{amsfonts}
\usepackage{amsmath}
\usepackage{amsthm}
\usepackage{amssymb}
\usepackage{amscd}
\usepackage[dvips]{graphicx}

\newcounter{num}
\newcommand{\Id}{\mathop{\rm Id}\nolimits}

\newcommand{\din}{\mathop{\rm deg_{in}}\nolimits}
\newcommand{\dout}{\mathop{\rm deg_{out}}\nolimits}
\newcommand{\FAT}[1]{\mbox{{$\mathbb{#1}$}}}

\newcommand{\zz}{\FAT{Z}}
\newcommand{\cc}{\FAT{C}}
\newcommand{\CC}{\FAT{C}}
\newcommand{\QQ}{\FAT{Q}}
\newcommand{\qq}{\FAT{Q}}

\newtheorem{lem}{\sc Lemma}
\newtheorem{prop}[lem]{\sc Proposition}
\newtheorem{cor}[lem]{\sc Corollary}
\newtheorem{thm}[lem]{\sc Theorem}
\newtheorem{df}[lem]{\sc Definition}

\input{amssym.def}
\input{amssym.tex}

\begin{document}

\title{A graph theoretic approach to graded identities for matrices}
\author{Darrell Haile\footnote{The first author was supported in part at the Technion by a fellowship from the Lady Davis Foundation.}\\
\small Department of Mathematics, Indiana University, Bloomington,
IN 47405 \\
\small haile@indiana.edu\\
\\ 
Michael Natapov %
\\
\small Department of Mathematics, Technion-Israel Institute of Technology, Haifa,
Israel 32000 \\
\small natapov@tx.technion.ac.il}

\date{}                                         

\maketitle


\begin{abstract}
We consider the algebra $M_k(\CC)$ of $k$-by-$k$ matrices over the complex numbers and view it as a crossed product with a group $G$ of order $k$ by imbedding $G$ in the symmetric group $S_k$ via the regular representation and imbedding $S_k$ in $M_k(\CC)$ in the usual way. This induces a natural $G-$grading on $M_k(\CC)$ which we call a crossed product grading. This grading is the so called elementary grading defined by any $k$-tuple $(g_1,g_2,\dots, g_k)$ of distinct elements $g_i \in G$. We study the graded polynomial identities for $M_k(\CC)$ equipped with a crossed product grading. To each multilinear monomial in the free graded algebra  we associate a directed labeled graph. This approach allows us to give new proofs of known results of Bahturin and Drensky on the generators of the $T$--ideal of identities and the Amitsur-Levitsky Theorem.

Our most substantial new result is the determination of the asymptotic formula for the $G$--graded codimension of $M_k(\cc)$.

\end{abstract}

\footnotetext{{\it 2000 Mathematics Subject Classification:} Primary 16S35, 16R10; Secondary 05C20.}

\footnotetext{{\it Keywords:} matrix algebra, graded polynomial identities, graded codimension, directed graph, Eulerian path.}


\section{Introduction.}

Let $G$ be an arbitrary group.  If $(g_1,g_2,\dots, g_k)$ is any $k$--tuple of elements of $G$ we can form an associated $G$--grading on $M_k(\cc)=\oplus_{g\in G}V_g$ as follows:  $V_g=0$ if $g$ is not of the form $g_i^{-1}g_j$ for some $i$ and $j$ and if $g=g_i^{-1}g_j$ for some $i$ and $j$, then $V_g$ is the span of the matrix units $e_{r,s}$ such that $g=g_r^{-1}g_s$. Such gradings are called elementary and have been studied by several authors, see e.g. \cite{BD}.  One is particularly interested in the graded polynomial identities for this grading,  polynomials $f(x_{i_1,h_1},\dots, x_{i_n,h_n})$ in variables indexed by elements of the group that become zero under all homogeneous substitutions, that is whenever each $x_{i_j,g_j}$ is replaced by an element of the component $V_{g_j}$.  For example in the case where  the $g_i$ are all distinct,  Bahturin and Drensky \cite{BD} have found an explicit set of generators for the $T$--ideal of all such identities.

In this paper we are interested in the case in which $G$ is a finite group, $k=|G|$ and the $g_i$ are distinct, so that the tuple $(g_1,g_2,\dots, g_k)$ is simply an ordering of the group elements.  It is shown in Aljadeff and Karasik \cite{AK}  that any two $k$--tuples of this type give $G$--isomorphic gradings and that the  grading comes from a crossed-product decomposition of $M_k(\cc)$.  We first recall the basic concepts.

Let $F$ be a field and $K$ a finite Galois extension of $F$ with Galois group $G$.  Given a two-cocycle $f:G\times G\rightarrow K^\times$ we can associate the crossed product algebra $K^fG=\oplus_{g\in G}Kx_g$, where the product is determined by the conditions $x_gk=g(k)x_g$ for all $g\in G$ and $k\in K$ and $x_gx_h=f(g,h)x_{gh}$ for all $g,h\in G$.  The algebra $K^fG$ is $F$--central simple of degree $k=|G|$ and every central simple $F$--algebra is Brauer equivalent to such a crossed product algebra for some choice of $G$, $K$ and $f$.  We will refer to such an algebra as a $G$--crossed product algebra.  Although $F$ is a field, $K$ need not be, and in fact we are particularly interested in the case where $F=\cc$, the field of complex numbers.  The only $\cc$--central simple algebra of degree $k$ is $M_k(\cc)$ and the only Galois extension of degree $k$ is $\cc^k$.  Any finite group $G$ of order $k$ may be viewed as a Galois group of $\cc^k$ over
$\cc$, where $G$ acts via the regular representation and the only cocycle is the trivial one, that is, $f(g,h)=1$ for all $g,h\in G$.  We can view this construction more concretely as follows:  Order the group elements as $e=g_1,g_2,\dots, g_k$ and label the matrix unit $e_{i,i}$ by $e_{g_i}$.
For each diagonal matrix $E = \sum_{i=1}^{k} a_i e_{g_i}$, let $E^g = \sum_{i=1}^{k} a_i e_{g g_i}$.
There is a homomorphism $\phi$ from the group $G$ into $P_k$, the group of permutation matrices in $M_k(\cc)$, that takes $g\in G$ to the permutation matrix $P_g=\phi(g)$ that satisfies $P_g E P_g^{-1}=E^g$, for all diagonal matrices $E$. We then have a decomposition $M_k(\cc)= \oplus_i D_kP_{g_i}$ where $D_k$ denotes the set of diagonal matrices.  The extension $D_k/\cc$ is our Galois extension and because $P_gP_h=P_{gh}$ for all $g,h\in G$, we have the desired crossed product structure (with trivial cocycle) on $M_k(\cc)$.
So $M_k(\cc)$ may be viewed as a $G$--crossed product for every group $G$ of order $k$.  Moreover Aljadeff and Karasik show in \cite{AK} that the matrix units $e_{i,j}$ are homogeneous and in fact $e_{i,j}\in D_kP_g$ if and only if $g=g_{i}^{-1}g_j$.  So this is precisely the grading determined by the tuple $(g_1,g_2,\dots,g_k)$.

We now return to graded identities.  We start with the free algebra $\qq\langle X^G\rangle$, where $X^G = \{x_{i,g} : 1\leq i, g\in G\}$. Each element $f$ in this algebra is a polynomial in the noncommuting variables $x_{i,g}$ with rational coefficients.  We evaluate such a polynomial on the $G$--crossed product $M_k(\cc)$, but allow only homogeneous evaluations.  In other words we can substitute for the variable $x_{i,g}$ elements from the component $D_kP_g$ only.  In particular we call a polynomial $f$ a {\it graded identity} for $M_k(\cc)$ if $f$ vanishes on every homogeneous substitution.  The set of these identities is an ideal in the free algebra and is an example of a {\it $T$--ideal}, which means that the ideal is stable under every graded endomorphism of the free algebra.

The main object of study in this paper are the {\it strongly multilinear} polynomials in the free algebra. Each such polynomial is a sum, with rational coefficients of monomials of the form $x_{i_1,g_1}x_{i_2,g_2}\cdots x_{i_n,g_n}$, where the subscripts $i_1,i_2,\dots , i_n$ are distinct.  We will refer to these as strongly multilinear monomials.
The adjective "strongly"  is to indicate that these monomials are not just multilinear in the variables $x_{i,g_i}$ but also in the numerical subscripts,  that is, we do not allow $x_{i,g}$ and $x_{i,h}$ to appear in the same monomial unless $g=h$. Using the process of linearization, it is easy to see that the $T$--ideal of graded identities is generated by the strongly multilinear graded identities it contains. Our main tool is a finite directed graph that we associate to each strongly multilinear monomial.  This graph has several interesting properties.  For example two strongly multilinear monomials have the same graph if and only if their difference is a graded identity. These differences, which we call binomial identities,  are basic to the theory of $G$--graded identities. In particular they generate the $T$--ideal of identities. We present new proofs of some known results to show the usefulness of the graph in Section \ref{section:graph}.

In Section \ref{section:asymptotics} we use the graphs to determine the asymptotic formula for the codimension of graded identities for $M_k(\cc)$.
 To put this in perspective, we first recall some results on the codimension growth in the ungraded case. Let $A$ be an algebra over a field $F$ of characteristic zero. Let $F\langle X \rangle$ be a free algebra on the set of countably many noncommuting variables $X = \{x_1, x_2, \dots, \}$, and let $\Id(A)$ denote the $T$--ideal of polynomial identities for $A$ in $F\langle X \rangle$. It is well known, that because $F$ is of characteristic zero,  $\Id(A)$ is completely determined by the multilinear identities. For each positive integer $n$ let $P_n$ be the vector space of multilinear polynomials of degree $n$:
 $$P_n={\rm span}\{x_{\sigma(1)}x_{\sigma(2)}\cdots x_{\sigma(n)} \  | \ \sigma \in S_n \},$$
 where $S_n$ is the symmetric group on the set $\{1,2,\dots , n\}$. The $n$-codimension $c_n(A)$ of the algebra $A$ is the dimension of $P_n$ modulo the identities:
 $$c_n(A) = \dim \frac{P_n}{P_n\cap \Id(A)}.$$
The codimension of a PI-algebra (i.e. an algebra satisfying a polynomial identity) was introduced by Regev in \cite{Regev1} where he proved that for any PI-algebra $A$ the codimension $c_n(A)$ is exponentially bounded. Regev conjectured in \cite{Regev3} that for any PI-algebra $A$ the asymptotic behavior of the codimension sequence $c_n(A)$ is given by
\begin{equation} \label{regev's congecture}
c_n(A) \sim a \cdot n^t \cdot \ell^n,
\end{equation}
where $a$, $t$ and  $\ell$ are some constants. Furthermore, in all cases computed so far, $\ell \in \zz$, $t \in \zz$, and $a \in \qq[\sqrt{2\pi}, \sqrt{b}]$ for some $0<b \in \zz$. In \cite{GZ1} and \cite{GZ2} Giambruno and Zaicev proved that for any PI-algebra $A$, the exponent of $A$, $\displaystyle {\rm exp}(A) = \lim_{n \rightarrow \infty} \sqrt[n]{c_n(A)}$ exists and is a nonnegative integer.  In \cite{BR} Berele and Regev proved that if $A$ is an algebra satisfying a Capelli identity (e.g. $A$ is finitely generated), then the conjecture (\ref{regev's congecture}) holds. In many cases an asymptotic formula for the codimension is known. In particular, Regev showed \cite{Regev2} that the codimension of $M_k(\CC)$ is given by

\begin{equation}\label{regev asymptotics}
 c_n(M_k(\CC)) \sim a \cdot n^{-\frac{k^2-1}{2}}\cdot k^{2n},
\end{equation}
where $\displaystyle a = \left(\frac{1}{\sqrt{2\pi}}\right)^{k-1} \left(\frac{1}{2}\right)^{\frac{1}{2}(k^2-1)}\cdot 1!2!\cdots (k-1)!k^{\frac{1}{2}(k^2+4)}.$

Now let $A$ be a an algebra graded by a group $G$. For each positive integer $n$ let $P^G_n$ be the vector space of strongly multilinear polynomials of degree $n$:
 $$P^G_n={\rm span}\{x_{\sigma(1), g_{\sigma(1)}}x_{\sigma(2),g_{\sigma(2)}}\cdots x_{\sigma(n),g_{\sigma(n)}} \  | \ \sigma \in S_n,\  g_1,g_2,\dots , g_n \in G \}.$$
The $G$--graded $n$-codimension $c^G_n(A)$ of the algebra $A$ is the dimension of $P^G_n$ modulo the graded identities:
 $$c^G_n(A) = \dim \frac{P^G_n}{P^G_n\cap \Id^G(A)}.$$

Much less is known about the graded codimensions. In \cite{AGM} Aljadeff, Giambruno and La Mattina proved that for a finite dimensional PI-algebra $A$ graded by an abelian group $G$ the graded exponent $\displaystyle {\rm exp}^G(A) = \lim_{n \rightarrow \infty} \sqrt[n]{c^G_n(A)}$ exists and is an integer. 
If $G$ is not abelian, then the exponent is not known in general. 
Aljadeff and Belov proved in \cite{AB} that if $A$ is an algebra with a fine $G$--grading (that is, $A$ is isomorphic to a twisted group algebra $F^cG$), then the the asymptotic behavior of $G$--graded codimensions of $A$ is given by $c^G_n(A) \sim |G'|\cdot |G|^n$, where $G'$ is the commutator subgroup of $G$. In particular, if the algebra $M_k(\CC)$ is fine graded by a group $G$ of order $k^2$, then its $G$--graded codimension is

\begin{equation}\label{aljadeff belov asymptotics}
 c^G_n(M_k(\cc)) \sim a \cdot k^{2n},
\end{equation}
where $a=|G'|$.

Now let $M_k(\cc)$ be equipped with a crossed product grading by a group $G$ of order $k$. We have the following asymptotic result:

\begin{equation}\label{haile natapov asymptotics}
 c^G_n(M_k(\cc)) \sim  a \cdot n^{-{k-1\over 2}} \cdot k^{2n},
\end{equation}
where $\displaystyle a = \left(\frac{1}{\sqrt{2\pi}}\right)^{k-1} \left(\frac{1}{2}\right)^{\frac{1}{2}(k-1)}k^{\frac{k}{2}+1}$.

One interesting thing about this formula is that, unlike the fine graded case, it depends only on the order of the group, that is, any two groups of the same order give the same asymptotics. In fact, we prove even more, namely, that the $G$--graded codimension itself does not depend on the group $G$. In other words it depends only on the algebra $M_k(\cc)$.

Another observation that arises from comparing the asymptotic formulas (\ref{regev asymptotics}), (\ref{aljadeff belov asymptotics}) and
(\ref{haile natapov asymptotics}) is that there is a certain pattern in the exponent of $n$. Notice that the exponent of $n$ in these formulas is a function of the dimension of the $e$-component of $M_k(\cc)$. Namely, we can regard the ungraded algebra as trivially graded by a degenerate group of order 1. In this case the homogeneous $e$-component is the whole algebra $M_k(\cc)$ of dimension $k^2$. In the fine graded case the $e$-component is of dimension 1, and in a crossed product grading the $e$-component is of dimension $k$. Thus, the exponent of $n$  is of the form $\displaystyle -{\dim M_k(\cc)_e-1\over 2}.$ It is natural to ask whether this is mere coincidence.


Our last results are two explicit formulas for the graded codimension for $M_2(\cc)$ derived using the graph count. One of these formulas (formula \ref{DV codimension}) was essentially established by Di Vincenzo in \cite{DV}. The other is a closed formula of the form:
$$\displaystyle c^{C_2}_n(M_2(\cc)) = \binom{2n+1}{n} - 2^n + 1.$$

\section{Binomial identities and the graph of a strongly multilinear monomial.}\label{section:graph}

We begin this section with a determination of a very useful set of generators for the $T$--ideal of $G$--graded identities for $M_k(\cc)$.   We first need a definition.

\begin{df} Let $m= x_{i_1,g_1}x_{i_2,g_2}\cdots x_{i_n,g_n}$ be a strongly multilinear monomial  in  $\qq\{x_{i,g} : 1\leq i, g\in G\}$.  For each $\pi\in S_n$ let $\pi(m)=x_{\pi(1),g_{\pi(1)}}x_{\pi(2),g_{\pi(2)}}\cdots x_{\pi(n)g_{\pi(n)}}$.  We will call $\pi$ an {\it initial product preserving} permutation for $m$ if

(1) $g_1g_2\cdots g_n=g_{\pi(1)}g_{\pi(2)}\cdots g_{\pi(n)}$, and

(2) For every $i$, $1\leq i\leq n$, if $\pi^{-1}(i)=j$, then $g_1g_2\cdots g_i=g_{\pi(1)}g_{\pi(2)}\cdots g_{\pi(j)}$.

\end{df}

In fact we will soon see that condition two implies condition one. We will call two monomials $m$ and $r$ in $\qq\{x_{i,g} : 1\leq i, g\in G\}$ {\it equivalent} if  $r=\pi(m)$ for some permutation $\pi$ that is initial product preserving for $m$. As in the introduction we let  $P^G_n={\rm span}\{x_{\sigma(1), g_{\sigma(1)}}x_{\sigma(2),g_{\sigma(2)}}\cdots x_{\sigma(n),g_{\sigma(n)}} \  | \ \sigma \in S_n,\  g_1,g_2,\dots , g_n \in G \}.$  Clearly the relation of equivalence is in fact an equivalence relation on the monomials in $P^G_n$.

\begin{prop}\label{prop:bimonial}  Let $f(x_{1,g_1}x_{2,g_2}\cdots x_{n,g_n})=\sum_{\pi\in S_n}a_\pi x_{\pi(1),g_{\pi(1)}}x_{\pi(2),g_{\pi(2)}}\cdots x_{\pi(n)g_{\pi(n)}} \in P^G_n$  be a strongly multilinear polynomial.   Let $f=f_1+f_2+\cdots +f_t$  be the decomposition of $f$ into the sums over equivalent monomials.   Then $f$ is a $G$--graded identity for $M_k(\cc)$ if and only if each $f_i$ is $G$--graded identity for $M_k(\cc)$. Moreover a given $f_i$ is an identity if and only if the sum of its coefficients is zero.
\end{prop}

\proof   Let $\displaystyle f(x_{1,g_1}x_{2,g_2}\cdots x_{n,g_n})=\sum_{\pi\in S_n}a_\pi x_{\pi(1),g_{\pi(1)}}x_{\pi(2),g_{\pi(2)}}\cdots x_{\pi(n)g_{\pi(n)}}$.    We evaluate $f$ on $M_k(\cc)$ by substituting, for each $x_{i,g_i}$,  the homogeneous element $E_iP_{g_i}$  where $E_i$ is a diagonal matrix.   We may choose the matrices $E_i$, $i=1,2,\dots, n$, with the property that the total set of $nk$  entries is algebraically independent over $\qq$.   Recall from the Introduction that if $E$ is diagonal and $g\in G$, then $P_gE=E^gP_g$.   Upon substitution we obtain $\displaystyle f(E_1P_{g_1},E_2P_{g_2}\dots E_nP_{g_n})= \sum_{\pi\in S_n}a_\pi E_{\pi(1)} E_{\pi(2)}^{g_{\pi(1)}}E_{\pi(3)}^{g_{\pi(1)}g_{\pi(2)}}\cdots E_{\pi(n)}^{g_{\pi(1)}g_{\pi(2)}\cdots g_{\pi(n-1)}}$.  Because the entries of the $E_i$ are  all algebraically independent  the only way this sum can be zero is if the  subsums over homogeneous elements with exactly the same resulting product of diagonal matrices is zero.  So the question is when two of these products, say
$E_{\pi(1)} E_{\pi(2)}^{g_{\pi(1)}}E_{\pi(3)}^{g_{\pi(1)}g_{\pi(2)}}\cdots E_{\pi(n)}^{g_{\pi(1)}g_{\pi(2)}\cdots g_{\pi(n-1)}}$  and $E_{\rho(1)} E_{\rho(2)}^{g_{\rho(1)}}E_{\rho(3)}^{g_{\rho(1)}g_{\rho(2)}}\cdots E_{\rho(n)}^{g_{\rho(1)}g_{\rho(2)}\cdots g_{\rho(n-1)}}$ (with $g_{\pi(1)}g_{\pi(2)}\cdots g_{\pi(n)}=g_{\rho(1)}g_{\rho(2)}\cdots g_{\rho(n)}$) are equal.   But the only way this can happen is if for each diagonal matrix $E_i$,  if  $E_i^\sigma$ appears in the first product and $E_i^\tau$   appears in the second product then $E_i^\sigma=E_i^\tau$.    It then follows that $\sigma$ must equal $\tau$.   This is exactly condition two of our definition of equivalence (where what we called $\pi$ there is now $\rho\pi^{-1}$) and so we are done.  \qed

In particular if $m$ and $\pi(m)$ are equivalent, then  $m-\pi(m)$ is a $G$--graded identity. We will call such identities {\it binomial} identities.

\begin{cor}  The space $P^G_n$ is spanned by the binomial identities  $m-\pi(m)$ for strongly multilinear monomials $m$ in $P^G_n$.   In particular the binomial identities generate the $T$--ideal of $G$--graded identities for $M_k(\cc)$.
\end{cor}

\begin{figure}[h]
$$\includegraphics[width=100mm]{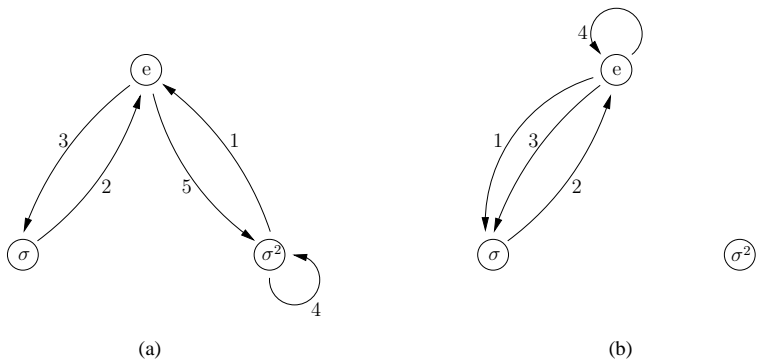}$$
  \caption{Let $G=\{e,\sigma, \sigma^2\}$ be the cyclic group of order 3, generated by $\sigma$. Given are (a) The graph of the monomial $x_{3,\sigma}x_{2,\sigma^2}x_{5,\sigma^2}x_{4,e}x_{1,\sigma}$. (b) The graph of the monomial $x_{3,\sigma}x_{2,\sigma^2}x_{4,e}x_{1,\sigma}$.}\label{fig:first2}
\end{figure}

We proceed to establish the connection of these identities to the theory of graphs. The kind of graphs we consider will be finite directed graphs, with labels on the vertices and the edges.  Every edge has a direction.  There may be several edges in both directions between two given vertices and there may be edges with the same beginning vertex and ending vertex.  So let $G$ be a finite group.  For each strongly multilinear monomial $x_{i_1,g_1}x_{i_2,g_2}\cdots x_{i_n,g_n}$ in $\qq\{x_{i,g} : 1\leq i, g\in G\}$ we construct a graph with vertices labeled by all of the elements of the group.  There is an edge labeled $i_1$ from the vertex labeled $e$ to the vertex $g_1$ and, for $j>1$, an edge labeled $i_j$ from the vertex labeled $g_1g_2\cdots g_{j-1}$ to the vertex $g_1g_2\cdots g_{j-1}g_j$.  In other words the graph is really a directed path through the vertices starting at $e$ and passing successively through $g_1, g_1g_2, g_1g_2g_3$, and so on, ending at $g_1g_2\cdots g_n$. We refer to such a path from $e$ to $g_1g_2\cdots g_n$ as an {\it Eulerian path}. The path may not hit all of the vertices and may hit the same vertex many times.  We exhibit some examples in Figure \ref{fig:first2}. 

It should be observed that every edge in the graph has an associated ``weight", the group element $g$ such that the initial vertex of the edge multiplied by $g$ gives the end vertex of that edge.  Of course many edges may have the same weight.  Moreover the group value of a vertex is equal to the product of the weights of the edges in any path from $e$ to that vertex.

From now on to simplify notation we will often write an ``arbitrary" strongly multilinear monomial as $x_{1,g_1}x_{2,g_2}\cdots x_{n,g_n}$ instead of $x_{i_1,g_1}x_{i_2,g_2}\cdots x_{i_n,g_n}$.

\begin{prop} \label{prop:the same graph} Let $m$ be a strongly multilinear monomial of degree $n$ and let $\pi\in S_n$.  Then $\pi(m)$ is equivalent to $m$ if and only if $m$ and $\pi(m)$ have the same graph.
\end{prop}

\proof Let $m=x_{1,g_1}x_{2,g_2}\cdots x_{n,g_n}$. Then $\pi(m)=x_{\pi(1),g_{\pi(1)}}x_{\pi(2),g_{\pi(2)}}\cdots x_{\pi(n)g_{\pi(n)}}$. The statement that $m$ and $\pi(m)$ have the same graph means that  in the graph for $m$, the sequence of edges labeled $\pi(1),\pi(2),\dots,\pi(n)$ is another path using each edge exactly once.  But as we observed above the group value of a vertex equals the product of the weights in any path in the graph from $e$ to that vertex.  Hence if $1\leq i\leq n$ and $\pi^{-1}(i)=j$ then when the new path determined by $\pi$ reaches $\pi(j)=i$ (that is, starts from $e$ and follows the edges $\pi(1),\pi(2),\dots \pi(j)$) we will have
$g_{\pi(1)}g_{\pi(2)}\cdots g_{\pi(j)}=g_1g_2\cdots g_i$, as desired.   This shows that if $\pi(m)$ is equivalent to $m$ then $m$ and $\pi(m)$ have the same graph.   It also shows that if $m$ and $\pi(m)$ have the same graph then condition two of the definition of equivalence is satisfied.   But the vertex $g_1g_2\cdots g_n$ is the last vertex in the path determined by $m$ and $g_{\pi(1)}g_{\pi(2)}\cdots g_{\pi(n)}$ is the last vertex in the path determined by $r=\pi(m)$.  This last vertex is uniquely determined by the property that it has one fewer edge leaving it than starting at it.  Hence if $m$ and $\pi(m)$ have the same graph both paths must end at the same point.  So $g_1g_2\cdots g_n=g_{\pi(1)}g_{\pi(2)}\cdots g_{\pi(n)}$  and $m$ is equivalent to $\pi(m)$. \qed

\begin{figure}[h]
$$\includegraphics[width=100mm]{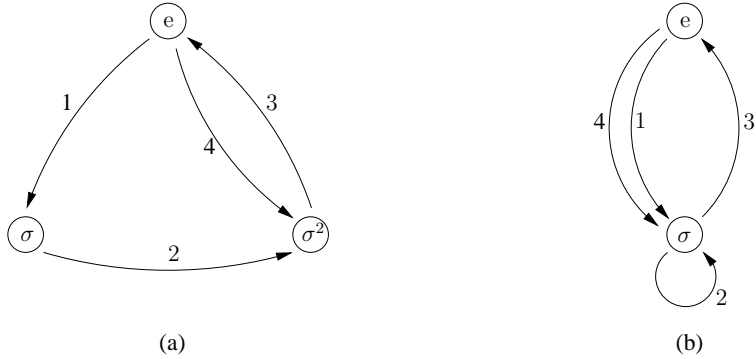}$$
  \caption{(a) The graph of the equivalent monomials $x_{1,\sigma}x_{2,\sigma}x_{3,\sigma}x_{4,\sigma^2}$ and $x_{4,\sigma^2}x_{3,\sigma}x_{1,\sigma}x_{2,\sigma}$, where $G$ is the cyclic group of order 3. (b) The graph of four equivalent monomials $x_{1,\sigma}x_{2,e}x_{3,\sigma}x_{4,\sigma}$, $x_{4,\sigma}x_{2,e}x_{3,\sigma}x_{1,\sigma}$, $x_{1,\sigma}x_{3,\sigma}x_{4,\sigma}x_{2,e}$, and $x_{4,\sigma}x_{3,\sigma}x_{1,\sigma}x_{2,e}$, where $G$ is the cyclic group of order 2.}\label{fig:first1}
\end{figure}

There are examples of graphs of equivalent monomials in Figure \ref{fig:first1}.\smallskip

\begin{cor}
In the definition of initial product preserving permutation,  condition two implies condition one.
\end{cor}

 We have the following obvious corollary:

\begin{cor}\label{cor:the same graph}
Let $m=x_{1,g_1}x_{2,g_2}\cdots x_{n,g_n}$ and $r=\pi(m)$.  The following conditions are equivalent:

(a) $m$ and $r$ are equivalent.

(b) $m$ and $r$ have the same graph.

(c) $m-r$ is a $G$--graded identity of $M_k(\cc)$.
\end{cor}

We want to investigate when two monomials are equivalent.  In other words we want to understand the initial product preserving permutations for a given monomial.  We first show that every such permutation is a product of very simple ones.  We will describe this simple kind of initial product preserving permutation in terms of the graph.  If in the graph of $m$ there are two vertices $g$ and $h$ (not necessarily distinct) for which there are two different segments of the path going from $g$ to $h$ then we can switch the order of these two segments.  This will give a new path and so this permutation will be initial product preserving.  We will call this a  {\it basic} permutation for $m$.
In Figure \ref{fig:first1}(a) we can switch the path segment from $e$ to $\sigma^2$ consisting of the edges 1 and 2 and the edge 4 to get the monomial  $x_{4,\sigma^2}x_{3,\sigma}x_{1,\sigma}x_{2,\sigma}$ equivalent to $x_{1,\sigma}x_{2,\sigma}x_{3,\sigma}x_{4,\sigma^2}$. In Figure \ref{fig:first1}(b) there are four (including the trivial one) basic permutations for the monomial $x_{1,\sigma}x_{2,e}x_{3,\sigma}x_{4,\sigma}$.

\begin{prop}
Let $m$ be a monomial.  Every initial product preserving permutation of $m$ is a product of basic permutations.
\end{prop}

 \proof Let $m=x_{1,g_1}x_{2,g_2}\cdots x_{n,g_n}$ and let $\pi$ be an initial product preserving permutation for $m$.  So $\pi$ determines another path (call it the $\pi$--path) through the graph of $m$.  We will show that we can find a basic permutation $\sigma$ such that the permutation $\sigma\pi$ satisfies $\sigma\pi(1)=1$.  The result will follow by induction on the degree of $m$.
 We may assume $\pi(1)\not= 1$.  The edge labeled 1 is one of edges in the $\pi$--path $\pi(1),\pi(2),\dots, \pi(n)$.
 First assume the edge labeled 2  comes before edge 1 in this path.  Because edge 2 begins at the endpoint of edge 1 this means that in this $\pi$--path, before we reach edge 1, there is a segment from $e$ to the endpoint of edge 1.  But when we do reach edge 1 in the $\pi$--path there is another segment from $e$ to the endpoint of edge 1, namely edge 1 itself.  We let $\sigma$ be the basic permutation that switches these two segments.  The composition $\sigma\pi$ then puts edge 1 back in the first position and finishes this case.
 So we may assume that edge 2 comes after edge 1 in the $\pi$--path.  Because edge 1 is not the first edge in the $\pi$--path there must then be an edge $r$, $r\geq 3$,  such that edge $r$ comes before edge 1, but edge $r-1$ comes after edge 1.  So there is a segment from $e$ to the initial point of edge $r$ (which is the same as the endpoint of edge $r-1$) and this segment comes before you reach edge 1.  Because edge $r-1$ comes after edge 1, there is then another segment from the initial point  of edge 1, that is the vertex $e$, to the endpoint of edge $r-1$.  So we have two path segments from $e$ to the endpoint of edge $r-1$ and we can let $\sigma$ be the basic permutation that switches these two segments.  The composition $\sigma\pi$ then puts edge 1 back in the first position and we are done.\qed

 \begin{cor}{\rm (Bahturin and Drensky)}
 The $T$--ideal of $G$--graded identities of $M_k(\cc)$ is generated by the following set of identities:

 (1) $x_{1,e}x_{2,e}-x_{2,e}x_{1,e}$

 (2) For each $g\in G$, $x_{1,g}x_{2,g^{-1}}x_{3,g}-x_{3,g}x_{2,g^{-1}}x_{1,g}$.

\end{cor}

\proof  It suffices to show that if $m$ is a strongly multilinear monomial and $\pi$ is an initial product preserving permutation of $m$ then $m-\pi(m)$ is in the $T$--ideal generated by elements of type (1) and (2).
By the proposition we may assume $\pi$ is a basic permutation.  So there are vertices $g$ and $h$ such that there are two path segments from $g$ to $h$ and $\pi$ is the permutation that switches these two segments. Assume first that $g=h$.  In that case the two path segments are loops beginning and ending at $g$.  The product of the weights of the edges (in order) around each loop equals e.  Therefore on the monomial $m$, switching the two loops has the effect of switching two successive partial products both of which equal e.  Hence it is a consequence of an identity of type (1).  If $g\not = h$, then we have a path segment from $g$ to $h$ followed by a path segment from $h$ to $g$ followed by another path segment from $g$ to $h$.  The product of the weights of the edges (in order) from  $g$ to $h$ is $g^{-1}h$  while the product of the weights of the edges from $h$ to $g$ is $h^{-1}g$.  Hence in $m$ we have three successive partial products $g^{-1}h$, $h^{-1}g$, $g^{-1}h$ and the effect of  $\pi$ is to switch the first and third of these three segments.  But this is a consequence of an identity of type (2). \qed

Bahturin and Drensky actually prove a more general result.  We digress briefly to explain how our methods can be used to prove their more general statement.  They consider an arbitrary (not necessarily finite) group $G$ and the elementary grading on $A=M_k(\cc)$  that comes from a $k$--tuple $(g_1,g_2,\dots, g_k)$ of distinct elements of $G$. One really deals only with the subgroup of $G$ generated by the $g_i's$ and so we may assume $G$ is finite or countable.  We consider first the finite case, say $|G|=m$.    In that case we can extend the $k$--tuple to an $m$--tuple $(g_1,g_2,\dots, g_m)$ including all of the elements of $G$.  We then form the crossed-product grading determined by this $m$--tuple.  If we decompose $M_m(\cc)= \oplus_{i=1}^m D_mP_{g_i}$ as described in the Introduction, then we may identify $A$ with $eM_m(\cc)e$ where $e=e_{11}+e_{22}+\cdots +e_{kk}$.  In other words $A=\oplus_{i=1}^m eD_mP_{g_i}e=\oplus_{i=1}^mee^{g_i}D_mP_{g_i}$.  Notice that in this sum the nonzero terms are those of the form $ee^gD_mP_g$, where $g=g_i^{-1}g_j$ for some $i,j$, $1\leq i,j\leq k$.  We are interested in the $G$--graded identities on $A$.  As before the $T$--ideal of graded identities is generated by the strongly multilinear identities.  The new phenomenon here is that we may have monomial identities.  For example $x_{r,g}$ is an identity if $g$ is not of the form $g_i^{-1}g_j$ for some $i,j$, $1\leq i,j\leq k$.  If we have a strongly multilinear identity $f$ in which no monomial is an identity then the analysis we used before Proposition \ref{prop:the same graph} shows that $f$ is a linear combination of differences $r-\pi(r)$ where $r$ is a strongly multilinear monomial and $\pi$ is an initial product preserving permutation of $r$.  In particular these identities are generated by the basic identities of Corollary 6.  Notice that such an identity is also a $G$-graded identity for the crossed-product grading on $M_m(\cc)$.  So we are left with considering monomial identities.  Because our result seems to be more precise than that of Bahturin and Drensky, we will state it more formally:

\begin{prop}
Let $r$ be a strongly multilinear monomial that is an identity for $A=M_k(\cc)$.  Then $r$ is of the form $r=st$ where $t$ is an arbitrary monomial and $s$ is obtained by the $T$--operation from a monomial identity of degree at most $k$.
\end{prop}
\proof Let $r=x_{1,h_1}x_{2,h_2}\cdots x_{u,h_u}$ be a strongly multilinear monomial.  Then $r$ is an identity for $A$ if and only if $0=eD_mP_{h_1}D_mP_{h_2}\cdots D_mP_{h_u}e=ee^{h_1}e^{h_1h_2}\cdots e^{h_1h_2\cdots h_u}D_mP_{h_1h_2\cdots h_u}$.  Hence $r$ is an identity if and only if $ee^{h_1}e^{h_1h_2}\cdots e^{h_1h_2\cdots h_u}=0$.  Now let $E=\{g_1,g_2,\dots, g_k\}$.  If $g\in G$ let $Eg=\{\sigma g|\sigma\in E\}$.  The condition that $ee^{h_1}e^{h_1h_2}\cdots e^{h_1h_2\cdots h_u}=0$ is equivalent to the condition that $E\cap Eh_1\cap Eh_1h_2\cap \cdots \cap Eh_1h_2\cdots h_u=\emptyset$.  If for some $j$, $E\cap Eh_1\cap Eh_1h_2\cap \cdots \cap Eh_1h_2\cdots h_{j-1}=E\cap Eh_1\cap Eh_1h_2\cap \cdots \cap Eh_1h_2\cdots h_{j-1}\cap Eh_1h_2\cdots h_j$, then we may remove the term $Eh_1h_2\cdots h_j$ and still have an empty intersection.  Continuing in this way we obtain $E\cap Eh_1h_2\cdots h_{j_1}\cap Eh_1h_2\cdots h_{j_2}\cap \cdots \cap Eh_1h_2\cdots h_{j_n}=\emptyset$, with $n\leq |E|=k$.  Hence if we let $y_i=h_1h_2\cdots h_{j_i}$ for $1\leq i\leq n$, the monomial $x_{1,y_1}x_{2,y_2}\cdots x_{n,y_n}$ is an identity of degree at most $k$. Moreover the monomial $s=x_{1,h_1}x_{2,h_2}\cdots x_{n,h_n}$ is obtained from $x_{1,y_1}x_{2,y_2}\cdots x_{n,y_n}$ by the $T$--operation (and so is an identity) and $r=st$ where $t=x_{n+1,h_{n+1}1}x_{n+2,h_{n+2}}\cdots x_{u,h_u}$,
so we are done.  \qed

The case where $G$ is (countably) infinite can be treated in the same way by using $G$ to produce a crossed-product grading on $M_{\infty}(\cc)$, the algebra of column-finite matrices.  The rest of the discussion applies with only cosmetic changes. So we obtain the full Bahturin-Drensky result \cite[Theorem 4.5]{BD}:  The $T$ ideal of identities of $A$ is generated by the the identities of Corollary 6, where $g\in G$ is chosen with nonzero component $A_g$, and by the (finitely many) monomial identities of degree at most $k$.

 Not every monomial has a non-identity initial product preserving permutation.  For example in the cyclic group of order 3 generated by $\sigma$ the monomial $x_{1,\sigma}x_{2,e}x_{3,\sigma}x_{4,e}x_{5,\sigma}$ has no such permutation, as is easily checked.  However if the degree of the monomial $m$ is at least $2k$ where $k$ is the order of the group, then we can prove there is always a nontrivial initial product preserving permutation for $m$.  We include a proof because it is easy and shows the usefulness of the graph.  But in fact we will see soon that considerably more is true.

 \begin{prop}
 If $G$ has order $k$ and $m$ is a strongly multilinear monomial of degree at least $2k$ then there is a nontrivial initial product preserving permutation for $m$.
 \end{prop}

 \proof Because the graph of $m$ has at least $2k$ edges, either every  vertex is the endpoint of at least two edges or some  vertex is the endpoint of at least three edges.  In the first case the identity $e$ is reached twice.  Because the path begins at $e$ it follows that there will be two loops at $e$ and there is a basic permutation switching these two loops.  If some vertex is the endpoint of at least three edges, then there    will be two loops that both begin and end at that vertex and so there is a basic permutation switching these two loops.\qed

 Here is the real theorem.

 \begin{thm}\label{swan.thm} If $G$ has order $k$ and $m$ is a strongly multilinear monomial of degree at least $2k$ then there are an even number of initial product preserving permutation for $m$, half of them odd and half of them even.
 \end{thm}

 \proof In \cite{Sw}, \cite{Sw2} Swan proved that in any finite directed graph with $k$ vertices, if you are given a path of length $n$, $n\geq 2k$, from vertex $a$ to vertex $b$ and you label the successive edges in that path by $1,2,\dots,n$ then the number of permutations $\pi\in S_n$ such that $\pi(1),\pi(2),\dots, \pi(n)$ is another path from $a$ to $b$ is even and half of the permutations are odd and half of them are even. If we apply this to the graph of the monomial $m$ we immediately obtain the result. \qed

 We return to the theory of identities.  Let $f$ be a (nongraded) homogeneous multilinear polynomial over $\qq$ of degree $k$.  Such a polynomial is of the form $f(x_1,x_2,\dots,x_n)=\sum_{\pi\in S_k}a(\pi)x_{\pi(1)}x_{\pi(2)}\cdots x_{\pi(n)}$, where the coefficients $a(\pi)$ are rational numbers.  Because $f$ is multilinear, if we want to check whether it is an identity for $M_k(\cc)$ it suffices to check it on a basis of $M_k(\cc)$.  In particular if $G$ is a group of order $k$ and we look at the crossed-product decomposition $M_k(\cc)= \oplus D_kP_g$ then in order to show $f$ is an identity it suffices to evaluate it on homogeneous elements.  That is it suffices to show $f(t_1P_{g_1},t_2P_{g_2},\dots, t_nP_{g_n})=0$ where the $t_i's$ are arbitrary elements in $D_k$ and the $g_i's$ are arbitrary elements of $G$.  But this is the same as saying that for every choice of $g_i's$, $f(x_{1,g_1},x_{2,g_2},\dots, x_{n,g_n})$ is a $G$--graded identity for $M_k(\cc)$. By Proposition \ref{prop:bimonial}, if we partition the monomials in this expression using the equivalence relation we see that to be an identity we must have that the sum of the $a_{\pi}'s$ in each class equals zero.

 For example take $f$ to be the standard polynomial $\displaystyle s_n= \sum_{\pi\in S_k}sgn(\pi)x_{\pi(1)}x_{\pi(2)}\cdots x_{\pi(n)}.$  The analysis above shows that $s_n$ is an identity for $M_k(\cc)$ if and only if $s_n(x_{1,g_1},x_{2,g_2},\dots, x_{n,g_n})$ is a $G$--graded identity for every choice of $g_i's$ in $G$ and this is true if and only if each monomial appearing in $s_n(x_{1,g_1},x_{2,g_2},\dots, x_{n,g_n})$ has an even number of initial product preserving permutations, half of them odd and half of them even. But if $n\geq 2k$ then this statement is true by Theorem \ref{swan.thm}.  Hence if $n\geq 2k$, then $s_n$ is an identity for $M_k(\cc)$.  This is the Amitsur-Levitsky theorem.  However our proof is only partly new.  The main ingredient is Swan's theorem, which Swan proved  precisely to give a new proof of Amitsur-Levitsky.  His use of his graph theorem however did not involve graded identities.

 It should be pointed out that there is a purely group theoretic formulation of the previous theorem:  Let $G$ be a group of order $k$ and let $g_1g_2\cdots g_n$ be a word in elements of $G$.  We define an initial product preserving permutation for the word as above.  That is, a permutation $\pi\in S_n$ is an initial product preserving permutation for the word $g_1g_2\cdots g_n$ if for every $i$, $1\leq i\leq n$, if $\pi^{-1}(i)=j$, then $g_1g_2\cdots g_i=g_{\pi(1)}g_{\pi(2)}\cdots g_{\pi(j)}$.  Then the statement is that if $n$ is at least $2k$ then there will be an even number of initial product preserving permutations for $g_1g_2\cdots g_n$, half odd and half even.  Moreover, this statement is equivalent to the Amitsur-Levitsky theorem.

\section{The graded codimensions.}\label{section:asymptotics}

In this section we analyze the asymptotic behavior of the codimension growth of the graded identities of $M_k(\CC)$.



Let $G$ be a group of order $k$. Let ${\bf g} = {\bf g}(k,n)$ denote a directed graph on $k$ vertices labeled by the elements of the group $G = \{e=g_1, g_2, \ldots, g_k\}$ and with $n$ edges labeled by the positive integers $\{1,2, \ldots, n\}$.

Recall that an {\it Eulerian path} from the vertex $g_i$ to the vertex $g_j$ is the enumeration $i_1, i_2, \ldots, i_n$ of all the edges of ${\bf g}(k,n)$ such that the first edge $i_1$ starts at $g_i$, the last edge $i_n$ ends at $g_j$, and for all $1 \leq l \leq n-1$ the initial point (vertex) of the edge $i_{l+1}$ is the endpoint of the edge $i_l$. We refer to such a path as a {\it cycle} if $g_i = g_j$. Note that we do not require an Eulerian path or cycle to hit all the vertices of the graph ${\bf g}$. We denote $M_k(n)$ the set of all graphs ${\bf g}(k,n)$ which have an Eulerian path from the vertex $e$ to the vertex $g_i$, for some $1\leq i \leq k$, and denote $|M_k(n)| = m_k(n)$.

Consider the $G$--crossed product grading on $A = M_k(\CC)$. Recall that the $G$--graded $n$-codimension $c^G_n(A)$ of the algebra $A$ is the dimension of the space $P_n^G/(P_n^G\cap \Id^G(A))$. Here $P_n^G$ is the $\qq-$vector space spanned by the monomials
$x_{i_1,g_1}x_{i_2,g_2}\dots x_{i_n,g_n}$, where $\{i_1,i_2, \ldots,i_n\} = \{1,2, \ldots,n\}$ and $g_1, g_2, \ldots, g_n$ are arbitrary
elements of $G$, and $\Id^G(A)$ is the $T$-ideal of graded polynomial
identities for $A$.

We first establish one-to-one correspondence between the equivalence classes of monomials in $P_n^G/(P_n^G\cap \Id^G(A))$ and the directed labeled graphs ${\bf g}(k,n) \in M_k(n)$. This will prove the following theorem:

\begin{thm}\label{independent}
Let $G$ be a group of order $k \geq 2$.
Then the graded $n$-codimension $c^G_n(A) = m_k(n)$. In particular, $c^G_n(A)$ does not depend on the group $G$.
\end{thm}

\proof We have seen in Section \ref{section:graph} that  a strongly multilinear monomial $x_{i_1,g_1}x_{i_2,g_2}\cdots x_{i_n,g_n}$ in $\QQ\langle X^G \rangle$ gives rise to a graph ${\bf g}(k,n)$ on $k$ vertices labeled by the elements of the group $G$ and with $n$ edges labeled by the integers $\{1,2, \ldots, n\}$. Moreover the edges $i_1, i_2, \dots, i_n$ constitute an Eulerian path starting at $e$ in this graph. Clearly, the converse is also true, namely, any graph ${\bf g}(k,n)$ with an Eulerian path from the vertex $e$ to the vertex $g_i$, for some  $1\leq i \leq k$, corresponds to a strongly multilinear monomial in $\QQ\langle X^G \rangle$. Moreover, we have seen in Proposition \ref{prop:the same graph} that two strongly multilinear monomials are equivalent modulo the $T$--ideal $\Id^G(A)$ if and only if they represent two Eulerian paths from the vertex $e$ (to the same vertex $g_i$) in the same graph ${\bf g}(k,n)$. In other words, the $G$--graded $n$-codimension $c^G_n(A)$ of $A$ is equal to the number of different graphs ${\bf g}(k,n)$ with an Eulerian path starting at $e$.\qed

Because of this theorem we denote the $G$--graded $n$-codimension $c^G_n(A)$ by $c_k(n)$.

In order to compute the $n$-codimention $c_k(n)$ we need to count the number of graphs in the set $M_k(n)$. In practice, we can give an exact formula for $m_2(n)$, that is, in case $k=2$ only. In the general case our strategy is to count the number of graphs of a more general kind and to show that the number of these graphs has the same asymptotic behavior as $m_k(n)$.

We now describes the more general kind of graphs we need. Given a graph ${\bf g}(k,n)$ as above, we say that the {\it degree in} of the vertex $g_i$, $\din(g_i)$, is the number of edges with the terminal point $g_i$, and the {\it degree out} of $g_i$, $\dout(g_i)$, is the number of edges with the initial point $g_i$. It is well known that a graph $\bf g$ has an Eulerian path from $g_i$ to $g_j$ if and only if

1. $\bf g$ is either connected, or the union of a connected subgraph and isolated points $g_l$ of $\din(g_l) = \dout(g_l) = 0$, where $l \neq i,j$.

2. for all $l \neq i,j$,  $\din(g_l)=\dout(g_l)$,

3. if $i=j$, then $\din(g_i)=\dout(g_i)$, and

4. if $i \neq j$, then $\din(g_i)=\dout(g_i)-1$ and $\din(g_j)=\dout(g_j)+1$.

\begin{figure}[h]
$$\includegraphics[width=120mm]{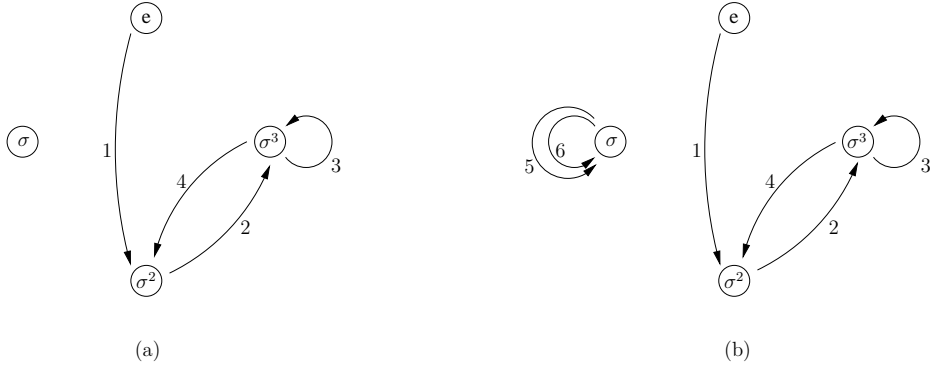}$$
  \caption{(a) Weakly connected graph with an Eulerian path from $e$ to $\sigma^2$. It is the graph of the monomial $x_{1,\sigma^2}x_{2,\sigma}x_{3,e}x_{4,\sigma^3} \in \QQ\langle X^G\rangle$, where $G$ is the cyclic group of order 4. (b) Strongly disconnected graph with an Eulerian pseudo-path from $e$ to $\sigma^2$. It does not correspond to any monomial in $\QQ\langle X^G\rangle$.} \label{fig:def}
\end{figure}

We say that a graph is {\it weakly connected} if it satisfies condition 1, and that a graph is {\it strongly disconnected} if it is not weakly connected. We say that a graph $\bf g$ has an Eulerian {\it pseudo-path} from $g_i$ to $g_j$ if it satisfies conditions 2-4, but is not necessarily weakly connected. We say that a graph $\bf g$ is {\it balanced} if it has an Eulerian pseudo-cycle, that is $\din(g_i)=\dout(g_i)$ for all $1\leq i\leq k$, but is not necessarily weakly connected. See Figure \ref{fig:def} for an illustration of these notions.\smallskip

We introduce the following notation.
Let $\Gamma_k(n)$ be the set of all graphs ${\bf g}(k,n)$ which have an Eulerian pseudo-path from the vertex $e$ to the vertex $g_i$, for some $1\leq i \leq k$, and denote $|\Gamma_k(n)| = \gamma_k(n)$. Let $P_k(n)$ be the set of all balanced graphs ${\bf g}(k,n)$, and denote $|P_k(n)| = p_k(n)$.


We first establish the following connection between the number of graphs in the sets $\Gamma_k(n)$ and $P_k(n)$:

\begin{lem}\label{balanced-unbalanced}
$$\gamma_k(n) = \frac{1}{k}p_k(n+1).$$
\end{lem}

\proof Let $\widetilde{P}_k(n+1)$ be the set of all balanced graphs ${\bf g}(k,n+1)$ such that the endpoint of the edge labeled by $n+1$ is $e$. Since in a balanced graph with $n+1$ edges the edge $n+1$ may have any terminal point $g_i, 1\leq i\leq k,$ with equal probability, we have $\displaystyle |\widetilde{P}_k(n+1)| = \frac{1}{k}|P_k(n+1)|.$

We build a 1-1 correspondence between the sets $\widetilde{P}_k(n+1)$ and $\Gamma_k(n)$. Given a graph ${\bf g}(k,n+1)$ in $\widetilde{P}_k(n+1)$ we just erase the edge $n+1$. Then we get a graph ${\bf g}(k,n)$ in $\Gamma_k(n)$. This correspondence is 1-1 since any graph ${\bf g}(k,n)$ having an Eulerian pseudo-path from the vertex $e$ to the vertex $g_i$, for some  $1\leq i \leq k$, can be completed in a unique way to a balanced graph ${\bf g}(k,n+1) \in \widetilde{P}_k(n+1)$ by adding the edge $n+1$ from $g_i$ to $e$. Thus we have
$$|\Gamma_k(n)| = |\widetilde{P}_k(n+1)| = \frac{1}{k}|P_k(n+1)|,$$
and the lemma follows. \qed

Next we count the balanced graphs on $k$ vertices with $n$ edges in two different ways:

\begin{lem}
The number of balanced labeled graphs ${\bf g}(k,n)$ is given by
\begin{equation}\label{balanced}
  \displaystyle p_k(n) =
\sum_{n_1=0}^{n} \binom{n}{n_1}^2\sum_{n_2=0}^{n - n_1}\binom{n-n_1}{n_2}^2 \dots \sum_{n_{k-1}=0}^{n - n_1 - \ldots - n_{k-2}}
\binom{n - n_1 - \ldots - n_{k-2}}{n_{k-1}}^2 = \end{equation}
\begin{equation}\label{multinomial}
\sum_{n_1 + n_2 + \ldots + n_k = n}  \binom{n}{n_1\ n_2\ \dots \ n_k}^2.
\end{equation}
\end{lem}

\proof Given a set of $n$ edges labeled by the numbers $\{1,2, \ldots, n\}$ one defines a unique balanced graph by the following sequence of choices:

1. Choose $n_1$ edges having the initial point $g_1$, and, independently, choose $n_1$ edges having the terminal point $g_1$ to ensure $\dout(g_1) = \din(g_1)$. This can be done in $\displaystyle \binom{n}{n_1}^2$ ways.

2. Out of the $n-n_1$ edges that do not start at $g_1$ choose $n_2$ edges having the initial point $g_2$. Independently, out of the $n-n_1$ edges that do not end at $g_1$ choose $n_2$ edges having the terminal point $g_2$.

3. Act similarly for $g_3, \ldots, g_{k-1}$.

4. Note that the edges having $g_k$ as the initial or the terminal point are uniquely defined by the previous steps.

This proves the equation (\ref{balanced}).

Alternatively, for any partition $n_1 + n_2 + \ldots + n_k = n$ one can choose $n_1, n_2, \dots, n_k$ edges having the initial point $g_1, g_2, \dots, g_k$, respectively, and, independently, choose $n_1, n_2, \dots, n_k$ edges having the terminal point $g_1, g_2, \dots, g_k$, respectively. This will give the equation (\ref{multinomial}). \qed

Note that the expression (\ref{balanced}) can be written as
\begin{equation}\label{recursive 1}
p_k(n) = \sum_{i=0}^{n} \binom{n}{i}^2 p_{k-1}(i).
\end{equation}
More generally, in \cite{RS} Richmond and Shallit show that $p_{k_1+k_2}(n)$ can be written in terms of  $p_{k_1}$ and $p_{k_2}$:
\begin{equation}\label{recursive 2}
p_{k_1+k_2}(n) = \sum_{i=0}^{n} \binom{n}{i}^2 p_{k_1}(i)p_{k_2}(n-i).
\end{equation}

In addition, Richmond and Shallit obtain in \cite[Theorem 4]{RS} the asymptotic behavior of $p_k(n)$ using the sum (\ref{multinomial}):

\begin{thm}\label{Richmond} Let $k$ be an integer $\geq 2$. Then, as $n \rightarrow \infty$, we have
$$\displaystyle p_k(n) = \sum_{n_1 + n_2 + \ldots + n_k = n}  \binom{n}{n_1\ n_2\ \dots \ n_k}^2 \sim  k^{2n+\frac{k}{2}} (4\pi n)^{\frac{1-k}{2}}.$$
\end{thm}

We now show that the numbers $m_k(n)$ and $\gamma_k(n)$ have the same asymptotics:

\begin{prop}\label{the same asymptotics}
$$ \frac{m_k(n)}{\gamma_k(n)} \rightarrow 1, \ \mbox{as} \ n \rightarrow \infty.$$
\end{prop}

\proof To prove the proposition we show that the number of the strongly disconnected graphs in $\Gamma_k(n)$ becomes negligible as $n \rightarrow \infty$.  Namely, let $sd_k(n)$  denote the number of strongly disconnected graphs on $k$ vertices with $n$ edges with an Eulerian pseudo-path starting at the vertex $e$, and let $sc_k(n)$ denote the number of connected graphs on $k$ vertices with $n$ edges with an Eulerian path starting at $e$. We claim that
$$ \frac{sd_k(n)}{\gamma_k(n)} \rightarrow 0, \ \mbox{as} \ n \rightarrow \infty.$$
It will follow then
$$\displaystyle \frac{m_k(n)}{\gamma_k(n)}  = \frac{\gamma_k(n) - sd_k(n)}{\gamma_k(n)}\rightarrow 1.$$

We first find $sd_k(n)$. Let $j$, $1 \leq j \leq k-1$, be the number of vertices in the connected component of the vertex $e$. There are  $\displaystyle \binom{k-1}{j-1}$ ways to choose the non-$e$ vertices of the $e$-component.  One needs at least $j-1$ edges to ensure that a graph on $j$ vertices is connected. Moreover, the number of edges in the $e$-component cannot be $n$, since otherwise the graph is weakly connected. Let $i$, $j-1 \leq i \leq n-1$, be the number of edges in the connected $e$-component, and choose $i$ edges out of the $n$ possible. Then the number of different configurations of the $e$-component is $sc_{j}(i)$. The remaining part of a graph with a connected $e$-component as above is a balanced graph on $k-j$ vertices with $n-i$ edges. Thus
$$ sd_k(n) = \sum_{j=1}^{k-1} \binom{k-1}{j-1} \sum_{i=j-1}^{n-1} \binom{n}{i} sc_{j}(i)p_{k-j}(n-i).$$

We now fix a subset $\{e=g_1, g_{i_2}, \ldots, g_{i_j}\}$ of the group $G$. Let $sd^j_k(n)$ be the number of strongly disconnected graphs with a connected $e$-component on the vertices $\{e=g_1, g_{i_2}, \ldots, g_{i_j}\}$. Then
          $$sd^j_k(n) = \sum_{i=j-1}^{n-1} \binom{n}{i} sc_{j}(i) p_{k-j}(n-i).$$
Since $sc_{j}(i) \leq \gamma_{j}(i)$ we have:
          $$sd^j_k(n) \leq \sum_{i=0}^{n-1} \binom{n}{i} \gamma_{j}(i) p_{k-j}(n-i).$$
By Lemma \ref{balanced-unbalanced}, we have
          $$sd^j_k(n) \leq \sum_{i=0}^{n-1} \binom{n}{i} \frac{1}{j} p_{j}(i+1) p_{k-j}(n-i) = $$
          $$\sum_{i=0}^{n-1} \frac{i+1}{n+1} \binom{n+1}{i+1} \frac{1}{j} p_{j}(i+1) p_{k-j}(n-i).$$
We rewrite the last expression using $l=i+1$:
          $$sd^j_k(n) \leq \sum_{l=1}^{n} \frac{l}{n+1} \binom{n+1}{l} \frac{1}{j} p_{j}(l) p_{k-j}(n+1-l).$$
Since $\displaystyle \frac{l}{j(n+1)} \leq 1$ for all $1 \leq l \leq n+1$ , we have
          $$sd^j_k(n) \leq \sum_{l=1}^{n} \binom{n+1}{l} p_{j}(l) p_{k-j}(n+1-l).$$
Now, since $\displaystyle \binom{n+1}{l} \geq n+1$ for all $1 \leq l \leq n$, we may write
          $$sd^j_k(n) \leq \sum_{l=1}^{n} \frac{1}{n+1} \binom{n+1}{l}^2 p_{j}(l) p_{k-j}(n+1-l) \leq $$
          $$\frac{1}{n+1} p_{k-j}(n+1) + \sum_{l=1}^{n} \frac{1}{n+1} \binom{n+1}{l}^2 p_{j}(l) p_{k-j}(n+1-l) + \frac{1}{n+1} p_{j}(n+1) = $$
          $$\frac{1}{n+1}\sum_{l=0}^{n+1} \binom{n+1}{l}^2 p_{j}(l) p_{k-j}(n+1-l).$$
Thus, by equation (\ref{recursive 2}) we have
          $$sd^j_k(n) \leq \frac{1}{n+1}p_{k}(n+1).$$
It follows that
          $$\frac{sd^j_k(n)}{\gamma_k(n)} \leq \frac{\frac{1}{n+1}p_{k}(n+1)}{\frac{1}{k}p_k(n+1)} = \frac{k}{n+1} \rightarrow 0,$$
as $n \rightarrow \infty$.
Hence we have
          $$\displaystyle \frac{sd_k(n)}{\gamma_k(n)} = \sum_{j=1}^{k-1} \binom{k-1}{j-1} \frac{sd^j_k(n)}{\gamma_k(n)} \rightarrow 0.$$

This completes the proof of the proposition. \qed

\bigskip

We are now ready to prove the main theorem of this section:

\begin{thm}\label{asymptotic} Let $G$ be a group of order $k \geq 2$. Then, as $n \rightarrow \infty$, the $G$--graded $n$-codimension
of $M_k(\CC)$ equipped with the $G$--crossed product grading is
$$\displaystyle c_k(n) \sim \frac{ k^{\frac{k}{2}+1} }{(4\pi)^{\frac{k-1}{2}}}\ n^{-\frac{k-1}{2}}\ k^{2n}.$$
\end{thm}

\proof As we have mentioned above $c_k(n) = m_k(n)$. Hence, by Proposition \ref{the same asymptotics}, $c_k(n) \sim \gamma_k(n)$. By Lemma \ref{balanced-unbalanced}, $\displaystyle c_k(n) \sim \frac{1}{k}p_k(n+1)$. By Theorem \ref{Richmond},

$$\displaystyle c_k(n) \sim \frac{1}{k}\left(k^{2(n+1)+\frac{k}{2}} (4\pi (n+1))^{\frac{1-k}{2}}\right) = \frac{ k^{\frac{k}{2}+1} }{(4\pi)^{\frac{k-1}{2}}}\ (n+1)^{-\frac{k-1}{2}}\ k^{2n},$$
and the theorem follows.\qed

In case $k=2$ we are able to give an explicit value of the $n$-codimension:

\begin{thm}\label{asymptotic-case-2} Let $C_2$ be a cyclic group of order $2$. Then the $C_2$--graded $n$-codimension of $M_2(\CC)$
equipped with the $C_2$--crossed product grading is
        $$\displaystyle c_2(n) = \binom{2n+1}{n} - 2^n + 1.$$
In particular, as $n \rightarrow \infty$, we have
        $$\displaystyle c_2(n) \sim \frac{1}{\sqrt{\pi}} n^{-\frac{1}{2}}\ 2^{2n+1}.$$
\end{thm}

\proof Let $C_2 = \langle e, \sigma \rangle$ be a cyclic group of order 2. The number of balanced graphs on the vertices $e$ and $\sigma$ with $n$ labeled edges is
$\displaystyle \sum_{i=0}^{n} \binom{n}{i}^2$ (for each possible $0\leq \dout(e) \leq n$, choose $i = \dout(e)$ edges starting at $e$, and, independently, choose $i = \din(e) = \dout(e)$ edges ending at $e$). Among these graphs, exactly $2^n$ are disconnected. Among the disconnected graphs one (that having all the edges starting and ending at $e$) is weakly connected. The number of graphs on the vertices $e$ and $\sigma$ with $n$ labeled edges with an Eulerian path from $e$ to $\sigma$ is $\displaystyle \sum_{i=1}^{n} \binom{n}{i}\binom{n}{i-1}$ (for each possible $1\leq i \leq n$, choose $i = \dout(e)$ edges starting at $e$, and, independently, choose $i-1 = \din(e) = \dout(e)-1$ edges ending at $e$). Thus we have
        $$c_2(n) = \sum_{i=0}^{n} \binom{n}{i}^2 + \sum_{i=1}^{n} \binom{n}{i} \binom{n}{i-1} - 2^n +1.$$
Direct computation shows that
        $$\sum_{i=0}^{n} \binom{n}{i}^2 + \sum_{i=1}^{n} \binom{n}{i} \binom{n}{i-1} = \sum_{i=0}^{n}
          \binom{n}{i} \binom{n+1}{i+1} = \sum_{i=0}^{n}\binom{n}{i} \binom{n+1}{n-i}.$$
A standard argument shows that the last sum
        $$\sum_{i=0}^{n}\binom{n}{i} \binom{n+1}{n-i}  = \binom{2n+1}{n},$$
and the first statement of the theorem follows.

Now,  since $\displaystyle \binom{2n+1}{n} = \frac{1}{2}\binom{2n+2}{n+1}$, using Stirling's formula $\displaystyle n! \sim \sqrt{2\pi n} \left(\frac{n}{e}\right)^n$, as $n \rightarrow \infty$, one easily gets the second statement of the theorem.\qed

An additional expression for $c_2(n)$ arises from a different approach to counting the number of graphs. Namely, we first assign the weights $e$ and $\sigma$ to the edges and then distribute the edges between the vertices $e$ and $\sigma$ according to their weights. Given $n$ edges,  we choose $m$ edges to have the weight $\sigma$ and the remaining $m-n$ edges to have the weight $e$. For each $0\leq m \leq n$ there are $\displaystyle \binom{n}{m}$ possible assignments of weights. Now, for each assignment, if $m \neq 0$, we choose exactly $\frac{m}{2}$, if $m$ is even, or $\displaystyle \lfloor \frac{m}{2} \rfloor+1$, if $m$ is odd, edges of weight $\sigma$ to go from the vertex $e$ to the vertex $\sigma$ to ensure the existence of an Eulerian path starting at $e$. In both cases the remaining $\displaystyle \lfloor \frac{m}{2} \rfloor$ edges of weight $\sigma$ go from $\sigma$ to $e$, and so there are $\displaystyle \binom{m}{\lfloor \frac{m}{2} \rfloor}$ possible choices. Independently, we distribute the $n-m$ loops of weight $e$ between the two vertices. There are $2^{n-m}$ ways to do so. We get
\begin{equation}\label{DV codimension} 2^{n-m}\binom{m}{\lfloor \frac{m}{2} \rfloor} \end{equation} different graphs. It should be mentioned that the number (\ref{DV codimension}) of non-equivalent monomials in variables $y_1, \dots, y_{n-m}$ of weight $e$ and variables $z_1, \dots, z_m$ of weight $\sigma$ was obtained by Di Vincenzo in \cite[Lemma 2]{DV}. If $m=0$, that is, if all the edges have weight $e$, there is exactly one (that having all the edges starting and ending at $e$) graph having an Eulerian path starting at $e$. Thus we have
$$c_2(n) = 1 +  \sum_{m=1}^{n} 2^{n-m} \binom{n}{m} \binom{m}{\lfloor \frac{m}{2} \rfloor}.$$

\end{document}